\documentclass[12pt]{article}
\usepackage{amssymb,latexsym}

 \newcommand{\ZZ}{\mathbb{Z}}
 \newcommand{\NN}{\mathbb{N}}

\def\auts{cellular automata}
\newcommand{\autb}{{cellular automaton }}
 \newcommand{\aut}{{cellular automaton}}
 
\newcommand{\az}{{A^\ZZ}}
\newtheorem{thm}{Theorem}[section]
\newtheorem{pro}{Proposition}[section]

\newtheorem{lem}{Lemma}[section]
\newtheorem{cor}{Corollary}[section]
\newtheorem{defi}{Definition}[section]
\newcommand{\saut}{\par \addvspace{\baselineskip}}

\newcounter{premarkc}

\begin{document}

\title {Cellular automata and Lyapunov exponents } 
\author{P. TISSEUR\\
 Institut de Math\'ematiques de Luminy \\
 UPR 9016 - 163, avenue de Luminy
Case 907\\
13288 Marseille Cedex 9 France }
\date{\ }
\maketitle

\begin{abstract}
The first definition of Lyapunov exponents (depending on a 
probability measure) for a one-dimensional cellular 
automaton were introduced by Shereshevsky in 1991.
The existence of an almost everywhere constant  value for each of the two exponents 
 (left and right), requires  
particular conditions for the  measure.
Shereshevsky establishes an inequality involving these two constants and the 
metric entropies of both the shift and the cellular automaton.
In this article we first prove that the two Shereshevsky's exponents exist for a 
more suitable class of measures,
then,  keeping the same conditions, we define new exponents, called average
Lyapunov exponents  smaller or equal to the first ones.
We obtain two inequalities: the first one is analogous to the Shereshevsky's  but concerns the average 
exponents; the second is the Shereshevsky inequality  but with more suitable assumptions.
These results are illustrated by two non-trivial examples, both proving that average exponents provide a better
 bound for the entropy, and one showing that the inequalities are  strict in general.
\end{abstract}

\section{Introduction}
A one-dimensional \autb (CA) denote by $F$ is a discrete mathematical idealization of a space-time
physical system. The space, called configuration space, consist of a discrete,   
 regular,  doubly infinite one-dimensional lattice with the property that 
each site can take a finite 
number of different values. A configuration is defined when every sites are fixed.  
 The discrete time is
represented by the action of a cellular automaton $F$ on this space. This action  
consist to change the value of a site considering only a finite number of  values situated 
in the neighborhood  of this site on the previous time step. We say that we apply a local rule.
The definition and the name  cellular automaton  were first given by Von Neumanmn and 
Ulam for modeling biological self reproduction.

For differential systems, the
Lyapunov exponents are essentially local properties and it is natural to
introduce a corresponding definition in the discrete frame of a cellular
automaton, defined by a local rule.

A perturbation in the configuration space is intuitively a change of values on 
some site. In \cite{Wo86},  Wolfram call perturbation a change of a finite number of site and 
study with computer the propagation of these changes. He call Lyapunov exponents the 
speed of these propagations and suspect that there exists 
 relations  between the spatial and temporal entropies and these exponents (see \cite{Wo86} pages 
261; 514).  The question of these relations  appears as one of the 20 general questions raised by Wolfram about 
cellular automata (see \cite{Wo86} page 172).
In \cite{Sh92} Shereshevsky gave a
mathematical definition of the Lyapunov exponents for a cellular automaton.
A left or right perturbation of a configuration become the set of all the configurations which 
 differ from the first one at  the right or left side of the central site.
The Shereshevsky definition of the Lyapunov exponents require to take the maximum speed of propagation on all the 
shifted configurations. 
Shereshevsky define the left and right Lyapunov exponents maps ($\lambda^+ , \lambda^-$) (see subsection 3.1)  which characterize the speed of propagation of these 
perturbations with respect to a 
 cellular automaton and shift-invariant measure.
 Requiring the $F$-ergodicity for the 
measure he  obtains that the maps have almost everywhere the same value and 
 note the two constants  $\lambda_\mu^+$ and $\lambda_\mu^-$.
Then if $\mu$ is also shift invariant, denoting by $h_\mu (F)$, (resp. $h_\mu (\sigma )$) the 
metric entropy of $F$ (resp. the metric entropy of the shift $\sigma $), Shereshevsky
establishes an inequality presumed by Wolfram and similar to the  Pesin one ( \cite{polli} or
\cite{Ru79}) in the differentiable case:
\begin{equation}\label{equa1}
h_\mu (F)\le h_\mu (\sigma ) (\lambda_\mu^+ +\lambda_\mu^-)
\end{equation}
where
$h_\mu (F)$ and  $h_\mu (\sigma )$ are respectively the metric entropy of the
cellular automaton $F$ and the metric entropy of the shift $\sigma$.
\medskip

{\it The main  
reason for continuing the work of Shereshevsky is that we know very few
examples of  cellular automaton with $F$-ergodic measure in general. The 
only exception are the expansive ones.}
\medskip

{\it Another reason  is that when a  cellular automaton has equicontinuous points
in the topological support of the measure, the measure can not be  $F$-ergodic.}

The {\bf Proposition 3.1} asserts  that these exponents $\lambda_\mu^+$ and $\lambda_\mu^-$ also exist in
the case of a shift-ergodic measure which is only $F$-invariant.
With these last conditions the uniform measure which is shift-ergodic is 
also invariant for 
every onto cellular automata.
More generally (see \cite{COPA74}) if $X$ is a mixing subshift of finite type and $F$ a cellular automaton such that $F(X)=X$ then the Parry measure on $X$ verify the new conditions. 
\smallskip

From Proposition 3.1, the new measure conditions implies that the exponents $\lambda_\mu^+$
and $\lambda_\mu^-$ only depend on the topological support $S(\mu )$ of  $\mu$.
To be precise they quantify the maximum  speed of the propagation 
of perturbation on the set $S(\mu )$.

We show by examples that for cellular automata with equicontinuous points
in $S(\mu )$ the exponents $\lambda^+_\mu$ and $\lambda^-_\mu$
are strictly positives under the new assumptions (see example 1),  
 although
 the metric entropie is  equal to $0$ (see 
{\bf Proposition \ref{Enul}}).
\smallskip

Next we define new Lyapunov exponents ($I^+_\mu$, $I_\mu^-$) called average Lyapunov exponents 
defined respected to an $F$-invariant and shift-ergodic measure.
 
From {\bf Proposition \ref{maxinf}}, we assume that the new exponents are smaller or equal to the first ones.

They are  equal to $0$ when exist 
equicontinuous points (see {\bf Proposition \ref{Enul}} and example 1).

Proving that the sum of these two exponents has a sense (see {\bf Proposition \ref{inde}}) we 
state 
{\it the main result named  {\bf Theorem 5.1} which gives the inequality: }
\begin{equation}
 h_\mu (F)\le h_\mu (\sigma )(I^+_\mu
+I_\mu^-). 
\end{equation} 
In  example 2 we show that the average Lyapunov exponents can be  strictly smaller 
than the Shereshevsky one, even if there is not equicontinuous point.

 Finally in {\bf Proposition 5.3} we establish a topological inequality :

If we denote
by $\mu_u$ the uniform measure on $A^\ZZ$ and by $h_{top}(F)$ the topological
entropy of an onto \autb $F$ : $A^\ZZ\to A^\ZZ$ we obtain 
\[
h_{top}(F)\le \log\#A(\lambda^+_{\mu_u}+\lambda^-_{\mu_u}). 
\]

We underline that when it is useful we only put one synthetic expression 
$\lambda_\mu^{\pm}$, $\Lambda_n^\pm$, $I^\pm _n$, etc.

\section{Preliminary }
\subsection{Symbolic systems and \auts}
Let $A$ be a finite set or alphabet. Denote by $A^*$ the set of all
concatenations of letters in $A$. These concatenations are called words. The
length of a word $u\in A^*$ is denoted by  $\vert u\vert$.
The set of bi-infinite sequences
$x=(x_i)_{i\in\ZZ}$ is denoted by $A^\ZZ$. A point $x\in A^\ZZ$ is called a
configuration. For $i\le j$ in $\ZZ$ we denote by $x(i,j)$ the word $x_i\ldots
x_j$ and by $x(p,\infty )$ the infinite sequence $(v_i)_{i\in\NN}$ such that for
all $i\in\NN$ one has $v_i=x_{p+i-1}$.
We endow $A^\ZZ$ with the product
topology. The shift $\sigma \colon A^\ZZ\to A^\ZZ$ is defined by :
$\sigma (x)=(x_{i+1})_{i\in \ZZ}$. For each integer $t$ and each word $u$, we
call cylinder the set $[u]_t=\{x\in \az : x_t=u_1\ldots ;x_{t+\vert
u\vert}=u_{\vert u\vert}\}$.
For this topology $A^\ZZ$ is a compact metric
space. A metric compatible with this topology can be defined by the distance
$d(x,y)=2^{-i}$ where $i=\min\{\vert j\vert \,\mbox{ such that } x(j)\ne y(j)\}$.
The dynamical system $(A^\ZZ ,\sigma )$ is called the full shift. A subshift $X$
is a closed shift-invariant subset $X$ of $A^\ZZ$ endowed with the shift
$\sigma$. It is possible to identify $(X,\sigma )$ with the set $X$. A language
$L$ is an arbitrary subset of $A^*$. Let $L_n$ be the set of words of length $n$
of $L$. The language associated to the subshift $X$ is $L(X)=\{u\in A^*\vert
\,\exists x\in X,\; x(i,i+\vert u\vert w-1)=u\}$. It is well known that $(X
,\sigma )$ is completely described by $L(X)$. If $\alpha =\{A_1,\ldots ,\,
A_n\}$ and $\beta =\{B_1,\ldots ,\, B_m\}$ are two partitions denote
by $\alpha \vee \beta$ the partition $\{A_i\cap B_j\, i=1\ldots n, \,\,
j=1,\ldots ,\, m\}$.

Consider a probability measure $\mu$ on the Borel sigma-algebra
$\cal{B}$ of $A^\ZZ$. If $\mu$ is $\sigma$-invariant then the topological
support of $\mu$ is a subshift denoted by $S(\mu )$. We denote by $\cal{M}(F)$
the set of all $F$-invariant probability measures and by $\# A$ the cardinal of
the set A. The uniform probability
measure on $A^\ZZ$ is the measure such that $\mu ([u]_t)=(\# A)^{-k}$ for all
integers $t$ and words $u\in A^k$.  The metric entropy $h_\mu (T)$ of a
transformation $T$
is an isomorphism invariant between two $\mu$-preserving
transformations; its definition can be found in \cite{Walt} and many  other
ergodic theory books. A \autb (CA) is a continuous self-map $F$ on
$A^\ZZ$ commuting with the shift. The Curtis-Hedlund-Lyndon theorem \cite{DGS75}
states  that for every cellular automaton $F$ there exist an integer $r$ and a
block map $f$ from $A^{2r+1}$ to $A$ such that: $F(x)_i=f(x_{i-r},\ldots ,x_i
,\ldots ,x_{i+r}).$
 The integer $r$ is called the radius of the cellular
automaton. If the block map of a cellular automaton is such that 
$F(x)_i=f(x_{i},\ldots ,\ldots ,x_{i+r})$, 
 the cellular automaton is called
one-sided and can be extended a map on a two-sided shift $A^\ZZ$ or a map on a
one-sided shift $A^\NN$. If $X$ is a subshift of  $A^\ZZ$ and one has
$F(X)\subset X$, the restriction of $F$ to $X$ determines a dynamical system
$(X,F)$; it is
called a \autb on $X$.
\section{Lyapunov exponents with shift-ergodic and $F$-invariant measure}
\subsection{The information propagation map}
Consider a \aut $(X,F)$ where $X$ is a subshift of $A^\ZZ$. Set
$W_s^+(x)=\{y\in X \vert\;\forall i\ge s ;\; y_i=x_i \}$ and $W_s^-(x)=\{y\in X
\vert ;\;i\le s\; y_i=x_i\}$. We claim that $W_s^+(x)$ is the set of
perturbations made by infinite blocks of points of $X$ located in the negative
coordinates of $x$. For any integer $n$ and $x$ in $X$ one has
\[
\tilde\Lambda^+_n(x)=\min\{s\ge 0 :\forall  1\le i\le n ,\, F^i(W_0^+(x))
 \subset W_s^+(F^i(x))\}, 
\]
\[
\tilde\Lambda^-_n(x)=\min\{s\ge 0 :\forall  1\le i\le n ,\, F^i(W_0^-(x))
\subset W_{-s}^-(F^i(x))\}.
\]
Then we define the two shift-invariant maps $\Lambda^\pm _n(x)=\max_{i\in\ZZ}\tilde{\Lambda}^\pm _n(\sigma^i(x))$.
\begin{premark}
Clearly $\tilde{\Lambda}_n^+$ and $\tilde{\Lambda}_n^-$ are two continuous functions bounded by $rn$. We have  changed a little bit the definition of Shereshevsky (see \cite{Sh92} pages 3) in order to clarify  
some proofs but this change does not affect the limits of the sequences 
$(\frac{\Lambda_n^\pm}{n})_{n\in\NN}$.
\end{premark}
\subsection{One proof of the existence of $\lambda^+_\mu$ and
$\lambda^-_\mu$ when $\mu$ is $\sigma$-ergodic.}
In this section we
prove that the limits of $(\frac{\Lambda_n^\pm}{n})_{n\in\NN}$  exist almost everywhere when 
$\mu$ is $\sigma$-ergodic  and
$F(S(\mu )\subset S(\mu )$(more suitable conditions) without using the subadditive ergodic 
theorem.
But with this new condition the maximum Lyapunov exponents are
rather topological than measure-theoretic quantities because they only depend on
the topological support $S(\mu)$. 
\begin{pro}\label{Lmax}
If $\mu$ is shift-ergodic and $F(S(\mu ))\subset S(\mu )$,
for
$\mu$-almost all $x$ in $X\supset S(\mu )$ the limits $\lim_{n\to\infty}\frac{
\Lambda^+_n(x)}{n}$ and $\lim_{n\to\infty}\frac{\Lambda^-_n(x)}{n}$ exist and
take constant values $\lambda^+_\mu$ and $\lambda^-_\mu$.
Moreover for almost all $x$ one has 
\[
\lambda_\mu^\pm =\lim_{n\to\infty}\frac{
\Lambda^\pm _n(x)}{n}=\lim_{n\to\infty}\max_{y\in S(\mu )}\frac{
\Lambda^\pm _n(y)}{n}=\max_{y\in S(\mu )}\limsup_{n\to\infty}
\frac{\Lambda^\pm _n(y)}{n}.
\]
\end{pro}
\begin{pproof}
We give only the proof for  $\lambda_\mu^+$.
For any $F$-invariant subshift $Y$, denote by $\hat{\Lambda}^+_n(Y)=\max_{x\in Y}\tilde{\Lambda}^+_n(x)$.

From \cite {Sh92} we have ${\tilde\Lambda}_{n+m}^+(x)\le{\tilde\Lambda}_n^+(x)+ \tilde{\Lambda}_m^+
(\sigma^{\tilde{\Lambda}_n^+(x)}\circ F^n(x))$ which implies that 
$(\hat{\Lambda}^+_n(Y))_{n\in\NN}$ is a subadditive sequence.
To finish the proof we need to show that for almost all $x$ we have  $\Lambda_n^+ (x)=\hat{\Lambda}^+_n(S(\mu ))$ which implies that for almost all $x$,
 the sequence  $(\frac{\Lambda^+(x)}{n})$ has a limit.
\smallskip

 We show that there exists a
set $G$ of full measure such that for any integer $n$, the map $\Lambda^+_n$ is
constant on $G$ and the value of this constant is $\hat{\Lambda}^+_n(S(\mu ))$.
Let $L(S(\mu ))$ be the language associated to $S(\mu )$ and let $u \in L(S(\mu
))$ be a word of length $2rn+r$. Clearly the map
$\tilde{\Lambda}^+_n$ is constant on the cylinder $[u]_0$. Put 
\[
V_n=\{u\in
L_{2rn+r}(S(\mu ))\mbox { such that } \tilde{\Lambda}^+_n([u]_0)
=\hat{\Lambda}^+_n(S(\mu ))\},
\]  
\[
\mbox{and }\,\, G_n=\{x\in S(\mu )\, \vert \, \exists i\in\NN \mbox{ such that } x(i,i+2rn+r)\in
V_n\}.
\]
It is easily seen that for any $x\in G_n$,
$\Lambda^+_n(x)=\hat{\Lambda}^+_n(S(\mu ))$. For any $n\in \NN$ the set $G_n$ is
$\sigma$-invariant. It contains a cylinder $[u]_{0}$ such that $u\in
V_n$, so $\mu (G_n)\ge \mu ([u]_0)>0$ and as $\mu$ is $\sigma$-ergodic,
 $\mu (G_n)=1$. Then $\mu(\cap_{n=1}^\infty G_n)=1$ and the map $\Lambda^+_n$
takes the value $\hat{\Lambda}^+_n(S(\mu ))$ on a set of full
measure. 
\end{pproof}
\subsection{Average Lyapunov exponents}
In this section we introduce the  average exponents $I^+_\mu$ and $I^-_\mu$ which represent an
average rate of propagation along the shift orbit for almost all points of
$X$.
We are going to show that these two exponents are less than or equal to their homologue
$\lambda^+_\mu$ and $\lambda^-_\mu$.



For any integer $n$, any point $x$, the map $I^-_n$ 
gives the minimum integer $m$ such that all the perturbations in the right side of 
$x_m$ never move until the central coordinate while the $n$ first iterations
; the exponent $I_n^+$ has a symmetric
definition. Formally,
\[
I^-_n(x)=\min\{ s\in\NN \, \vert\, \, \forall \, 1\le i\le n ,\,\vert \,
F^i(W^-_{s}(x))\subset W_0^-(F^i(x))\}, 
\]
\[
I^+_n(x)=\min\{ s\in\NN \,\vert \, \, \forall \, 1\le i\le n ,\,\vert\,
F^i(W^+_{-s}(x))\subset W_0^+(F^i(x))\}. 
\]
\vskip -1 true cm 

\begin{premark}
Clearly $I^+_n$ and $I^-_n$ are two continuous functions bounded by $rn$.
\saut
\end{premark}
Set $
I^+_{n;\mu}=\int_X I^+_n(x)d\mu (x)$ and $I^-_{n;\mu}=
\int_X I^-_n(x)d\mu(x).$ The Birkhoff's theorem implies that for almost all $x$ one has 
 $I^+_{n;\mu}=\lim_{n\to\infty}\sum_{i=-m}^m\frac{1}{2m+1}I^+_n(\sigma^i
(x))$ and $I^-_{n;\mu}=\lim_{n\to\infty}\sum_{i=-m}^m\frac{1}{2m+1}I^-_n(\sigma^i(x))$.
\begin{defi}
Call average Lyapunov exponents the limits
\[
I^+_\mu=\liminf_{n\to\infty}\frac{I^+_{n;\mu}}{n}\,\, \mbox{ and }\,\,
I^-_\mu=\liminf_{n\to\infty}\frac{I^-_{n;\mu}}{n}. 
\]
\end {defi}
\begin{pro}\label{maxinf}
If $\mu$ is $\sigma$-ergodic and $F(S(\mu ))\subset S(\mu )$, then
$I^+_\mu \le \lambda^+_\mu$ and $I^-_\mu \le \lambda^-_\mu $.
\end{pro}
\begin{pproof}
 By definition of $I^+_n (x)$, there exists $i\le n$ such that
$F^i(W_{-I^+_n(x)+1}^+(x))$

$\not\subset W_0^+(F^i(x))$. Hence for all $x$ we have
$\tilde{\Lambda}_n^{+}(\sigma^{-I_n^+(x)+1}(x))\ge I^+_n(x)-1$, then
$\Lambda^{+}_n(x)\ge I^+_n(x)-1$.  

We can write that 
$
\int_{S(\mu )}I^+_n(x)d\mu (x)\le \int_{S(\mu )}\left (
\Lambda^{+}_n(x)+1\right ) d\mu (x) 
$ which im\-plies that 
$$
\lim\inf_{n\to\infty}\frac{I^+_{n;\mu }}{n}\le
\lim\inf_{n\to\infty}\int_{S(\mu )}\frac{\Lambda^{+}_n(x)+1}{n}d\mu (x).
$$
Then using the dominated convergence
theorem we get 
\[
 I^+_\mu \le \int_{S(\mu
)}\lim_{n\to\infty}\frac{\Lambda^+_n(x)}{n}d\mu (x) =\lambda^+_\mu . 
\]
The proof is the same for $I^-_\mu$ and $\lambda^-_\mu$.
\end{pproof}
{\bf Question:}
We do not know examples of sequences
$(\frac{I^+_{n;\mu}}{n})_{n\in\NN}$ and $(\frac{I^-_{n;\mu}}{n})_{n\in\NN}$
which do not converge. Do they exist?
\section{Equicontinuity and Shannon-McMillan-\\Breiman theorem}
\begin{defi}
A cellular automaton has equicontinuous points (or Lyapunov stable points) if
and only if there exists a point $x$ in $X$ such that for all $\epsilon >0$,
there exists $\delta >0$, such that for all $y$ in $X$ with $d(x,y)<\delta$ then
$d(F^n(x),F^n(y))<\epsilon$ for any $n$. 
\end{defi}
\begin{defi}
Let $F$ be a cellular automaton with radius $r$. A word $B \in A^{2k+1}$ is
called blocking word if  for all $x$ in $X$ such that $x(-k,k)=B$, there exists
an infinite word sequence $v_n$, $\vert v_n\vert =2i+1\ge r$, such that
$F^n(x)(-i,i)=v_n$ for all  $n\in\ZZ^*$.
\end{defi}
\vskip -1 true cm 

\begin{premark}
If $B$ is a blocking word and if a point $x$ verifies $x(-k,k)=B$, then the
sequence $F^n(x)(-\infty,-i)$ does not depend on
$x(k,+\infty)$ because $2i+1\ge r$.
A blocking word completely disconnects the evolution of the
coordinates to its left and to its right.
This imply that a point with infinitively many occurrences of a blocking word 
is an equicontinuous point.
\end{premark}
The relation between equicontinuity
points and blocking words was established in \cite{Ku94}
 (see also  \cite{BT1}).
\begin{pro}\label{block}\cite{Ku94}\cite{BT1}
A \autb $F$ with radius $r$ acting on a transitive subshift $X$ has
equicontinuous points if and only if it has a blocking word.
\end{pro}
Let $\alpha$ be a finite partition of $X$, let $T$ be a measurable action on $X$
and $\mu$ be a $T$-invariant measure. Denote by $P_{n,\alpha }^{T}(x)$
the element of the partition $\alpha\vee T^{-1}\alpha\vee \ldots \vee
T^{-n}\alpha$ which contains $x$ and $h_\mu (T,\alpha)$ the metric entropy
of $T$ with respect to the partition $\alpha$.
\begin{thm}\label{shannon} (Shannon-McMillan-Breiman)
If $\mu$ is $T$-invariant, for almost all $x$ $\lim_{n\to\infty}\frac{-1}{n}
\log\mu (P_{n,\alpha}^T(x))$
exists and one has 
\[
\int_X\lim_{n\to\infty}\frac{-1}{n}\log\mu (P_{n,\alpha}^T(x))d\mu
(x)=h_\mu (T,\alpha ).
\]
If $\mu$ is a $T$-ergodic measure, then for almost all $x$
\[
\lim_{n\to\infty}\frac{-1}{n}\log\mu (P_{n,\alpha}^T(x))=h_\mu (T,\alpha ).
\]
\end{thm}
We give a new version of the Shannon-McMillan-Breiman theorem,
in the case of a one to one and onto action $T$, this new Proposition 
will be used in the proof of the main result.
\begin{pro}\label{vShannon}
Let $T$ be a one to one and onto action and  $\mu$ an ergodic measure.
If we denote by  $P_{n,m,\alpha}^T(x)$ the element of the  partition $\alpha\vee T\alpha
\ldots \vee T^m\alpha\vee T^{-1}\alpha \ldots \vee T^{-n}\alpha $ 
which contains $x$, we have  
\begin{equation}
\lim_{n+m\to\infty}\frac{-1}{n+m}\log\mu (P_{n,m,\alpha_p}^T(x))
=h_\mu (T,\alpha_p ). 
\end{equation}
\end{pro}
We only give a sketch of the proof, the complete proof appears 
in \cite{prepri}.

{\it Sketch of the proof:} (see \cite{pet} in which the similar 
proof for Theorem \ref{shannon} appears)  

Let $i(\alpha )(x)=-\log\mu (P^T_\alpha (x))$ and $i(\alpha /\beta )=
-\sum_{A\in\alpha}\log\mu(A/\beta )(x)\chi_A(x)$.
We use  the next definition for the metric entropy  
\[
h_\mu (\alpha ,T)=\lim_{n\to\infty}\int_X i(\alpha /\vee_{k=1}^nT^k\alpha )(x)
d\mu (x).
\]
Using the ``two sided'' version of the Birkhoff theorem 
$\lim_{n+m\to\infty}\sum_{k=-m}^n f\circ T^k (x)=\int_X f(x)d\mu (x)$
 with $f=\lim_{n\to\infty}i(\alpha /\vee_{k=1}^{n}T^{-k}\alpha )
$ and showing that 
$
P_{n,m,\alpha_p}^T(x)=i(\vee_{k=-n}^mT^k\alpha )(x)=\sum_{k=-m}^{n-1} 
i(\alpha /\vee_{j=1}^{n-k}T^{-j}\alpha )\circ T^k+i(T^{m-1}\alpha )
$
we obtain 
\[
\lim_{n+m\to\infty}i(\vee_{k=-n}^mT^k\alpha )(x)=\lim_{n\to\infty}\int_X i(\alpha /\vee_{k=1}^nT^k\alpha )(x)
d\mu (x). \hskip 2.5 true cm \Box
\]

\section{Main results}
The  proof of our principal result, Theorem \ref{principal} relies on
two propositions and two lemmas. Proposition \ref{inde} establishes that one can
treat independently the perturbations coming from the right and the
perturbations
coming from the left: this allows to sum the two exponents. Lemma
\ref{suitein} permits to split the general proof into two cases, Lemma \ref{Epo}
solves the first case and Proposition \ref{Enul} solves the second.
\begin{pro}\label{inde}
For any triple of positive integers $(n,p,i)$ with $i\le n$ and $p\ge r$
and for every $x$ in $X$, one has 
\[
F^i \left (C_{-p-I^+_n(\sigma^{-p}(x))}^{p+I^-_n(\sigma^{p}(x))}(x)\right )
\subset C_{-p}^p\left (F^i(x)\right ).
\]
\end{pro}
This means that if the point $y$ has the same coordinate as $x$
from $-p-I_n^+(\sigma^p(x))$ to $p+I^-_n(\sigma^{p}(x))$, then for each $i\le n$
the coordinates of $F^i(y)$ are equal of those
of $F^i(x)$ from $-p$ to $p$.

\begin{pproof}
Fix two positives integers $n$ and $p$. Choose a point $x\in X$ and
put $s^+ = I^+_n(\sigma^{-p}(x))$ and $s^- = I^-_n(\sigma^{p}(x))$. For each
point $y \in C^{p+s^-}_{-p-s^+}(x)$, set $y_1$ and $y_2$ such that
$y_1(-p-s^+,\infty )=x(-p-s^+,\infty )$, $y_1(-\infty ,p+s^-)=y(-\infty
,p+s^-)$, $y_2(-\infty ,p+s^- )=x(-\infty ,p+s^-)$ and $y_2(-p-s^+, \infty
)=y(-p-s^+, \infty)$. By definition of $I_n^+$ and $I_n^-$, for all $i\le
n$ one has  $F^i(y_1)(-p,\infty )=F^i(x)(-p, \infty )$ and $F^i(y_2)(-\infty ,
p)=F^i(x)(-\infty ,p)$. The proof consists in showing by induction that for
every
positive integer $i\le n$ one has 
\begin{equation}
F^i(y)(-p,p)=F^i(x)(-p,p).
\end{equation}
Recall that $f \colon A^{2r+1}\to A$ is the local map of
$F$; for every integer $k$ we also denote by $f$ the map from $A^{2r+1+k}$ to
$A^k$ define by
\[
 f(u_0\ldots u_{2r+k+1})=f(u_0\ldots u_{2r})f(u_1\ldots
u_{2r+1})\ldots f(u_{k-1}\ldots u_{2r+k}). 
\]
Let us prove the first step of the recurrence. If $p\ge r$ then by
$F(y)(-p,p)=F(x)(-p,p)$. As $y_1(-p-r,r)=y(-p-r,r)$ and $y_2(-r,p+r)=y(-r,p+r)$
one has 
\[
F(y)(-p,p)=f(y(-p-r,p+r))=f(y(-p-r,r))f(y(-r,p+r))
\]
\[
=f(y_1(-p-r,r))f(y_2(-r,p+r))
\]
using the definition of $y_1$ and $y_2$ we obtain
\[
F(y)(-p,p)=F(x)(-p,0)F(x)(0,p)=F(x)(-p,p).
\]
Let $i$ be a positive integer such that $i\le n-1$. We show that if (4) is
true for each $k\le i$ it remains true for $k=i+1$.

First we need the two equalities
\small
\begin{equation}
F^i(y_1)(-p-r,r)=F^i(y)(-p-r,r)\mbox{ and }F^i(y_2)(-r,p+r)=F^i(y)(-r,p+r).
\end{equation}
\normalsize
We give the proof of the first one, the second is analogous.
\saut

We prove the equality $F^i(y_1)(-p-r,r)=F^i(y)(-p-r,r)$ using a secondary
recurrence that establishes that for any positive integer $k\le i$ one has
\begin{equation}
 F^k(y_1)(-p-r(i+1-k),r)=F^k(y)(-p-r(i+1-k),r).
\end{equation}
As $y_1(-p-(r+1)i,p+s^-)=y(-p-(r+1)i,p+s^-)$, from the definition of $s^-$ one
has $F(y_1)(-p-ri,p)=F(y)(-p-ri,p)$ and since $p\ge r$ we obtain the first step
of this new recurrence. We suppose (6) is true for each $k\le i-1$, i.e.,
$F^k(y_1)(-p-r(i+1-k),r)=F^k(y)(-p-r(i+1-k),r)$; then
\[
F^{k+1}(y_1)(-p-r(i-k),r)=f(F^k(y_1)(-p-r(i+1-k), 2r))
\]
\[
=f(F^k(y_1)(-p-r(i+1-k), r)f(F^k(y_1)(-r,2r)). 
\]
and since (6) is true for each $k\le i-1$ we have 
$$
F^{k+1}(y_1)(-p-r(i-k),r)=f(F^k(y)(-p-r(i+1-k),r))f(F^k(y_1)(-r,2r))
$$
\[
=F^{k+1}(y(-p-r(i-k),0))f(F^k(y_1)(-r,2r)).
\]
To finish the proof of the step $k+1$ of this secondary recurrence, i.e.,
$F^{k+1}(y_1)(-p-r(i-k),r)=F^{k+1}k(y)(-p-r(i-k),r)$, we need to verify that
$f(F^k(y_1)(-r,2r))=f(F^k(y)(-r,2r))=F^{k+1}(y)(0,r)$.  From the definition of
$y_1$ we can assert that  $F^k(y_1)(-r,2r)=F^k(x)(-r,2r)$ and if we use the
hypothesis of the main recurrence (6), namely,
$F^{k+1}(y)(-p,p)=F^{k+1}(x)(-p,p)$ (since $k\le i-1$), we conclude that
\[
f(F^k(y_1)(-r,2r))=f(F^k(x)(-r,2r))=F^{k+1}(x)(0,r) 
\]
\[
 =F^{k+1}(y)(0,r).
\]
We are now in position to show that
$F^{i+1}(y)(-p,p)=F^{i+1}(x)(-p,p)$, which completes the proof of the main
recurrence. Starting from the equality
\small
$$
F^{i+1}(y)(-p,p)=f(F^i(y)(-p-r,p+r))=f(F^i(y)(-p-r,r))f(F^i(y)(-r,p+r)) 
$$
\normalsize
and using the two equalities (5) we conclude that 
\[
\begin{array}{ll}
F^{i+1}(y)(-p,p) \\
=f(F^i(y_1)(-p-r,r))f(F^i(y_2)(-r,p+r))\\
=f(F^i(x)(-p-r,r))f(F^i(x)(-r,p+r))\\
=F^{i+1}(x)(-p,0)F^{i+1}(x)(0,p)\\
=F^{i+1}(x)(-p,p).
\end{array}
\]
\end{pproof}
\begin {lem}\label {suitein}
Let $\mu$ be a shift-ergodic measure.
If there exist an equicontinuous point $x$ in $S(\mu )$ then 
for every integer $p\ge r$ and for almost all point the sequences 
$(I^+_n(\sigma^{-p}(x))+I_n^-(\sigma^p(x)))_{n\in\NN}$ are bounded.

If there is not equicontinuous point in $S(\mu )$ then for every 
$p\ge r$ and for almost all points $x$, the sequences $(I^+_n(\sigma^{-p}(x))+
I_n^-(\sigma^p(x)))_{n\in\NN}$ go to infinity.
\end{lem}
\begin{pproof}
If exist an equicontinuous point $x$ in $S(\mu )$ then there exist an integer 
$k$ and a  
blocking word $B=x(-k,k)$ (see Proposition \ref{block}).
Let $V(B)$ be the set of all the point with infinitely  many occurrences of
$B$ in the positive and negative coordinates.
From Remark 3 we claim that all the point of $V(B)$ are equicontinuous points.
As $\mu$ is shift-ergodic and $\mu ([B]_0)>0$ then  one has $\mu (V(B))=1$. 
Using Remark 3 we conclude that for each point $y\in V(B)$ and for 
each integer $p$ the sequences 
 $(I^+_n(\sigma^{-p}(y))+I_n^-(\sigma^p(y)))_{n\in\NN}$ are bounded.

We suppose now that there is not equicontinuous point in $S(\mu )$ and 
that there exist a set $E$ with 
strictly positive measure such that $E$ contains only points $x$  with the properties
 $\exists p(x)\in\NN \, \vert (I^+_n(\sigma^{-p}(x))+I_n^-(\sigma^p(x)))_{n\in\NN}$
is a bounded sequence.

Clearly $E\cap S(\mu )\neq \emptyset$. Let $y\in E\cap S(\mu )$. 
There exist $p\ge r$ such that  $M(+)=\max_{n\in\NN}\{I^+_n(\sigma^{-p}(y))\}$ and
$M(-)=\max_{n\in\NN}\{I^-_n(\sigma^p(y))\}$ are well defined.
From Proposition \ref{inde} for all $i\in\NN$ one has 
\[
F^i \left (C_{-p-M(-)}^{p+M(+)}(y)\right )
\subset C_{-p}^p\left (F^i(y)\right ).
\]
which implies that the word $B'=y_{-M(-)-p},\ldots
y_{M(+)+p}$ is a blocking word for $F$.
As $\mu$ is shift-ergodic and $\mu ([B']_0)>0$ then there exist a point 
$z\in S(\mu )$ with infinitely  many occurrences of 
$B'$ in the positive and negative coordinates. This point $z$ is an 
equicontinuous point (see Remark 3) which contradict the hypothesis.
\end{pproof}
\begin {lem}\label{Epo}
If $\mu$ is a shift-ergodic and $F$-invariant measure such that $F$ has no
equicontinuous point in $S(\mu )$, then 
$
\,\,h_\mu (F)\le h_\mu (\sigma ) (I^+_\mu +I^-_\mu). 
$
\end {lem}
\begin{pproof}
Fix $x\in X$ and denote by $\alpha_p$ the partition of $X$ into cylinders
$C_{-p}^p$ ($x\in X$); call $P_{n,\alpha_p}^F(x)$ the element of the
partition
$\alpha_p\vee F^{-1}\alpha_p\vee \ldots \vee F^{-n}\alpha_p$ that contains $x$.
By Lemma \ref{inde}, for any choice of positive integers $p$  $(p\ge r)$, $n$
and $i$ ($i\le n$), one has 
\[
 F^i \left
(C_{-p-I^+_n(\sigma^{-p}(x))}^{p+I^-_n(\sigma^{p}(x)})\right )\subset
C^p_{-p}\left (F^i(x)\right ).
\]
 The last inclusion implies that each $F^{-i}\alpha_p$ has an
element that contains the cylinder
$C_{-p-I^+_n(\sigma^{-p}(x))}^{p+I^-_n(\sigma^{p}(x))}(x)$, so 
$
P_{n,\alpha_p}^F(x)\supset C_{-p-I^+_n(\sigma^{-p}(x))}^{p+I^-_n(\sigma
^p(x))}(x),
$
 and consequently 
\begin{equation}
 -\frac{1}{n}\log \mu
(P_{n,\alpha_p}^F(x))\le -\frac{1}{n}\log \mu
\left(C_{-p-I^+_n(\sigma^{-p}(x))}^{p+I^-_n(\sigma^{p}(x))}(x)\right). 
\end{equation}
Applying the Shannon-McMillan-Breiman theorem to $F$ one shows that

$h_\mu (F, \alpha_p)=\int_X\lim_{n\to\infty}-\frac{1}{n}\log\mu (P_{n,\alpha_p}^F(x))$.
 Then by (7)
\[ 
h_\mu (F, \alpha_p)
\le \int_X \liminf_{n\to\infty} -\frac{1}{n}\log \mu \left
(C_{-p-I^+_n(\sigma^{-p}(x))}^{p+I^-_n(\sigma^{p}(x))}(x)\right )d\mu(x)
\,\,\,\mbox{ and}
\]
\small
$$
\hskip -3 true cm h_\mu (F, \alpha_p) 
\le \int_X \liminf_{n\to\infty}
-\frac{\log \mu
\left (C_{-p-I^+_n(\sigma^{-p}(x))}^{p+I^-_n(\sigma^{p}(x))}(x)\right)}{I^+_n(\sigma^{-p
}(x))+I^-_n(\sigma^{p}(x))}
$$
$$
\hskip 4 true cm \times\frac{I^+_n(\sigma^{-p}(x))+I^-_n(\sigma^{p}(x))}{n}
 d\mu (x).\eqno{(8)}
$$
\normalsize

By Lemma \ref {suitein} if there is no equicontinuous point in
$S(\mu )$ then for all integer $p$ and almost every point $x$ the sequence
$(I^+_n(\sigma^{-p}(x))+I^-_n(\sigma^p(x)))_{n\in\NN}$ goes to infinity.
 Considering that $\mu$ is
shift-ergodic and $\sigma$ is a one-to-one and onto map we can apply to $\sigma$ the version
(3) of the Shannon-McMillan-Breiman theorem, which gives
\[
 \lim_{n\to\infty} -\frac{\log \mu
\left(C_{-p-I^+_n(\sigma^{-p}(x))}^{p+I^-_n(\sigma^{p}(x)))}(x)\right)}{I^+_n(\sigma^{-p
}(x))+I^-_n(\sigma^{p}(x))} =h_\mu (\sigma ,\alpha_p)=h_\mu (\sigma )
\]
for almost all $x$ and every
positive integer $p$. Combining the last equality with (8)
 yields
\[h_\mu (F,
\alpha_p)\le h_\mu (\sigma )\times \int_X \liminf_{n\to\infty}
\frac{I^+_n(\sigma^{-p}(x))+I^-_n(\sigma^{p}(x))}{n}d\mu (x).
\] 
Using the Fatou lemma, we get 
\[
h_\mu (F, \alpha_p)\le
h_\mu (\sigma ) \times \liminf_{n\to\infty}\int_X\frac
{I^+_n(\sigma^{-p}(x))+I^-_n(\sigma^{p}(x))}{n}d\mu (x).
\]
Since $\alpha_p$ is an increasing sequence with the property $\bigvee_0^\infty
\alpha_i=\cal{B}$ and $\mu$ is $\sigma$-invariant we obtain
\[
\lim_{p\to\infty} h_\mu (F, \alpha_p)=h_\mu (F) \le h_\mu (\sigma )\times
\liminf_{n\to\infty}\int_X \frac{I^+_n(x)+I^-_n(x)}{n}d\mu (x).
\]
This last inequality completes the proof, so  
$\,\,h_\mu (F)\le h_\mu (\sigma )\times (I^+_\mu +I^-_\mu ).$
\end{pproof}
The next proposition establishes that if there exists a blocking word $u$ such
that $\mu ([u]_0)>0$ then the metric entropy $h_\mu(F)$ is equal to $0$.
\begin{pro}\label{Enul}
If a cellular automaton $F$ has equicontinuous points belonging to $S(\mu )$
then the average Lyapunov exponents $I^+_\mu$ and $I^-_\mu $ and the metric
entropy $h_\mu (F)$ are $0$. 
\end {pro} 
\begin{pproof} By our hypothesis and 
Lemma \ref{suitein}, for each integer $p\ge r$ the sequences
$(I^+_n(\sigma^{-p}(x))+I^-_n(\sigma^p(x)))_{n\in\NN}$ are bounded for almost all $x$. This
implies that
for any positive integer $p$ and for almost all $x$,
\[
\liminf_{n\to\infty}-\frac{\log\mu
\left (C_{-p-I^+_n(\sigma^{-p}(x))}^{p+I^-_n(\sigma^{p}(x))}(x)\right)}{I^+_n(\sigma^{-p
}(x))+I^-_n(\sigma^{p}(x))}
\] 
is bounded.
Then for all $p$ and for almost all $x$, we get
\[
\liminf_{n\to\infty}\frac {-\log\mu
\left(C_{-p-I^+_n(\sigma^{-p}(x))}^{p+I^-_n(\sigma^{p}(x))}(x)\right)}
{I^+_n(\sigma^{-p}(x))+I^-_n(\sigma^{p}(x))} \times
\frac{I^+_n(\sigma^{-p}(x))+I^-_n(\sigma^{p}(x))}{n}=0.
\]
From (8) in the proof of Lemma \ref{Epo}, the sum over $X$ of  the last equality is  an upper bound of the metric entropy $h_\mu (F)$ which implies that this
entropy is equal to $0$. On the other hand the sequence
$(\frac{I^+_n(x)}{n})_{n\in\NN}$ is bounded by $r$ and converges to $0$ for
almost all $x$, then applying the dominated convergence theorem one gets
\[
I^+_\mu =\lim_{n\to\infty} \int_X \frac{I^+_n(x)}{n}d\mu (x)\le
\int_X\lim_{n\to\infty} \frac{I^+_n(x)}{n}d\mu (x)=0. 
\]
 The proof is identical
for $I^-_\mu$. 
\end{pproof}
\vskip -1 true cm 

\begin{premark}
One can prove that $h_\mu (F)=0$ if there exist equicontinuous points in $S(\mu
)$ using Katok's definition of metric entropy.
\end{premark}

Combining Lemma \ref{Epo} and Proposition \ref{Enul} we obtain the next
theorem :
\begin {thm}\label{principal}
If $\mu$ is a $\sigma$-ergodic and $F$-invariant measure then
\[
h_\mu (F)\le h_\mu (\sigma ) (I^+_\mu +I^-_\mu ).
\]
\end {thm}
\vskip -0.7 true cm 
\begin{premark}
For one-sided cellular automata one defines a unique average
Lyapunov exponent $I_\mu$  whose definition
is identical to that of $I^+_\mu$ in this Subsection. Then the proof of the
inequality $h_\mu (F)\le h_\mu (\sigma )I_\mu$ does not require the use
of Proposition \ref{inde}.
\end{premark}
Since $\lambda_\mu^+\ge I^+_\mu$ and $\lambda_\mu^-\ge I^-_\mu$
(Lemma \ref{maxinf}) one has
\begin{cor}\label{corfi} If $\mu$ is a $\sigma$-ergodic and $F$-invariant measure then
\[
h_\mu (F)\le h_\mu (\sigma ) (\lambda_\mu^+ +\lambda^-_\mu ). 
\] 
\end{cor}
\subsection*{A topological inequality}
Here we recall some  definitions relative to the topological entropy that we 
denote by $h_{top}(F)$.  
Let $(X,F)$ be a dynamical system. 
 For any integer $n$ the distance $d_n$ is defined by $\forall x, y\in X^2$ one has $d_n(x,y)=\max\{d(F^i(x),F^i(y))\, 0\le i\le n\}$.
An $(n,\epsilon )$-covering set is a cover of $X$ by balls of diameter $\epsilon$ for the $d_n$ metric.
Let $D(n,\epsilon )$ be the minimum cardinal of an $(n,\epsilon$) covering set.
\[
h_{Top}(F)=\lim_{\epsilon \to\infty}\lim_{n\to\infty}\frac{1}{n}\log(D(n,\epsilon )).
\]
Let $\mu_u$ be the uniform measure of $A^\ZZ$. We will give a upper bound of
$h_{top}(F)$ according to the exponents $\lambda^+_{\mu_u}$ and
$\lambda^-_{\mu_u}$. Remark that for all cellular automaton $F$, the uniform measure  satisfies
the two conditions of Proposition \ref{Lmax}, so  $\lambda^+_{\mu_u}$ and
$\lambda^-_{\mu_u}$ always exist.
\begin{pro}\label{itop}
For any onto \autb $F$ : $A^\ZZ\to A^\ZZ$ one has 
$h_{top}(F)\le (\lambda^+_{\mu_u}+\lambda^-_{\mu_u})\log \# A.$
\end{pro}
\begin{pproof}
From Proposition \ref{inde} and proof of Proposition \ref{maxinf} ($\Lambda_n^\pm (x)+1\ge I^\pm _n(x)$), it follows that 
for any choice of positive integers $p$ $(p \ge r)$, $n$
and $i$ ($i\le n$), one has 
\[
 F^i \left
(C_{-p-\Lambda^{-}_n(x)-1}^{p+\Lambda_n^{+}(x)+1}(x)\right )\subset
C_{-p}^p (F^i(x)).
\]
Denote by  $\Omega (n,p))$ the set of all the cylinders $\left(C_{-p-\Lambda^{-}_n(x)-1}^{p+\Lambda_n^{+}(x)+1}(x)\,\, (x\in A^\ZZ)\right)$ and by $\hat{\Lambda}_n^{\pm}$ the maximum of all 
the $\Lambda_n^\pm (x)$.
The last inequality  implies that $\Omega (n,p))$ is a $(n,2^
{-p})$ covering set which show that  for all integers $n$ and $p$, we get $D(n,2^{-p})\le \# \Omega (n,p)$.
As 
\[
\#\Omega (n,p)=\#\{C_{-p-\hat{\Lambda}_n^{+}-1}^{p+\hat{\Lambda}_n^{-}+1}(x_j)\,
\vert \, x_j\in A^\ZZ\}
=(\#A)^{(2p+3+\hat{\Lambda}_n^{+}+\hat{\Lambda}_n^{-})},
\]
we can assert that  
$$
\hskip -4 true cm h_{Top}(F)=\lim_{p\to\infty}\lim_{n\to\infty}\frac{1}{n}\log D(n,2^{-p})
$$
$$
\hskip 3 true cm  \le 
\limsup_{p\to\infty}\lim_{n\to\infty}\frac{2p+3+\hat{\Lambda}_n^{+}
+\hat{\Lambda}_n^{-}}{n}\times \log \#A.
$$
Using Proposition \ref{Lmax}  we obtain 
$
\,\, h_{Top}(F)\le (\lambda^+_{\mu_u}+\lambda^-_{\mu_u})\log \#A.
$
 \end{pproof}

\section{Examples}
The two following examples show  that $I^+_\mu$ and $I^-_\mu$ can be strictly
less than $\lambda^+_\mu$ and $\lambda^-_\mu$.
The example 6.2 shows that  the inequality of 
Theorem \ref{principal} 
is in general strict.
In both examples we use the uniform measure which is shift-ergodic and
$F$-invariant when $F$ is onto from $A^\ZZ$ to itself.
\subsection{Coven's cellular automata}
In \cite{Co80}, Coven computes exactly the positive topological entropy of a
particular class of onto cellular automata with complex behavior. In
\cite{Bl96} Blanchard and Maass show that all these CA have equicontinuous
points.

A Coven aperiodic CA is defined by its block map $f : \{0,1\}^{r+1}\to
\{0,1\}$:

$
f(x_0,x_1,\ldots ,x_r)=(x_0+1)$ mod 2 if $x_1\ldots x_r=b_1\ldots b_r$,
$f(x_0,x_1,\ldots ,x_r)=x_0$ otherwise.
The word $B=b_1\ldots b_r$ must be aperiodic, which means that  for any integer $r>1$ there is no  integer $p$ ($0<p<r$) such that $b_{i+p}=b_i$ for $i=1,\ldots r-p$.
In \cite{Co80} Coven proves that the topological entropy of this type of CA is
$\log (2)$. Here we consider the Coven CA with radius $r=2$ and aperiodic
word $B=10$. This
particular example has the typical behavior of all the other 
Coven's automata. Let $\mu$ be the uniform measure
on $\{0,1\}^\ZZ$.  From
\cite{Bl96} we know that $000$ is a blocking word for $F$.  If
$\mu$ is the uniform measure $I^+_\mu +I^-_\mu=0$
by Proposition \ref{Enul} and $h_\mu(F)=0$.  On
the contrary the sum of the maximum Lyapunov exponents is strictly positive.
First it is clear that $\lambda^+_\mu =0$, because the block map $f$ does not
depends on  negative coordinates of $x$.

Let $y$ be the fixed point with $y_i=1$ for all $i$ and let $z$ be
a point with all the coordinates equals to $1$ except $z_0$. The word $01$ never
appears in $y$ so $F(y)=y$. Considering that $F(z)(-4,-2)=110$ we deduce that
$F^2(-6,-4)=110$ and by a trivial induction
$F^n(z)(-2n-2,-2n)=110$. Considering that $y$ belong to $S(\mu )=A^\ZZ$ and
applying Proposition \ref{Lmax} we see that $\lambda_\mu^-\ge\lambda^-(y)\ge 2$. The value
of $\lambda^-_\mu $ must be less than or equal to the radius of the (CA) which is equal to
$2$ so $\lambda^-_\mu=2$. It is well known that the topological entropy of the
two-shift is equal to $\log 2$. 
 From Theorem \ref{principal} we get $h_\mu (\sigma )(\lambda_\mu^+
+\lambda_\mu^-)=2\log 2>h_\mu (F)=0$. Remark that in this case the values of
$\lambda^+_\mu$ and $\lambda^-_\mu$ do not allow to prove that $h_\mu(F)=0$.
From \cite{Bl96} $h_{top}(F)=2\log (2)=h_\mu (\sigma )(\lambda_\mu^+
+\lambda_\mu^-)$ so in this case the inequality of Corollary \ref{corfi}
 becomes an equality.
\subsection{A sensitive cellular automaton}
\begin{defi}
Let $X$ a compact space and $T$ be a transformation of $X$. The
map $T$ is said to be sensitive if there exists a real $\epsilon >0$ such that
for any $x\in X$, any real $\delta >0$, there exists a positive integer $n$
and a
point $y$ such that $d(x,y)<\delta$ and $d(T^n(x),T^n(y))\ge \epsilon$.
\end{defi}
K\accent23urka \cite{Ku94} shows that a cellular automaton is sensitive if and
only if it has no equicontinuous points. By Proposition \ref{block} sensitive
cellular automata have no blocking words, so for all $x \in X$ one has
$\lim_{n\to\infty}(I_n^+(x)+I_n^-(x))=\infty$.

The aim of studying the sensitive cellular automaton $F$ defined below is
twofold. First, in spite of its rather simple behavior, it gives a good idea of
the reason why average Lyapunov exponents give a better upper bound of the metric
entropy $h_\mu (F)$. Secondly, this example shows that inequality (2) is
sometimes strict.
\saut

Set $X_1=\{0,1\}^\ZZ$, $X_2=\{0,1,2\}^\ZZ$ and $X=X_1\times
X_2$. Denote by $\mu_1$ the uniform measure on $X_1$, by $\mu_2$ the uniform
measure on $X_2$ and  $\mu$ the product measure $\mu_1\times \mu_2$ on $X$.
Clearly $\mu$ is the uniform measure on $X$, so $\mu$ is shift-ergodic. The
cellular automaton $F$ is the product of $F_1$ acting on $X_1$ and $F_2$ acting
on $X_2$. Denote by $\sigma$ the shift on $X$. The automaton $F_1$ is only
the shift on $X_1$.  For each $x \in X_1$ one has
$\tilde{\Lambda}^{F_1-}_n(x)=n$. As $\mu_1$ is shift-ergodic on $X_1$ and
$F_1$-invariant we can assert that $I^-_{\mu_1}=\lambda^-_{\mu_1}=1$.

We define a cellular automaton $F_2$ on $X_2$ with radius $r$ by its
local map $f_2$:
\[
f_2(x_{-r},\ldots x_0,\ldots x_r)=x_0+x_r\mbox{ if } 2\notin \{x_0,x_1,\ldots
,x_r\} 
\]
\[
\mbox{and }f_2(x_r,\ldots x_0,\ldots x_r)=x_0\mbox{ if } 2\in
\{x_0,x_1,\ldots ,x_r\}. 
\]
Using a criterion given in \cite{He69}, one can easily show that $F_2$ is onto, which implies that the
product automaton  $F$ is also onto.
As the uniform measure is invariant for an onto cellular automaton the exponents
$(I^+_\mu ,I^-_\mu )$ and $(\lambda^+_\mu ,\lambda^-_\mu )$
are defined for $F$ and $\mu$. Remark that
for each couple of integers $k$ and $i$ the value of $F^k(x)_i$ does not depend
on the coordinates at the left of $x_i$, so for  $F$ one has $I^+_\mu
=\lambda^+_\mu =0$. The letter $2$ is clearly a blocking word for
$F_2$. Every point $x\in X_2$ with infinitely many occurrences of $2$ in the
negative and positive coordinates is an equicontinuous point for $(F_2,X_2)$.
The restriction of $F_2$ to the subshift $\{0,1\}^\ZZ$ is the $r$ times iterated
shift. It follows that for $\mu_2$-almost all $x$ one has
$\Lambda^{-F_2}_n(x)=rn$, and by Proposition \ref{Lmax} $\lambda^-_{\mu_2}=r$. The
measure $\mu_2$ which is shift-ergodic on $X_2$ is also $F_2$-invariant because
$F_2$ is an onto map from $X_2$ to $X_2$. As $F_2$ has equicontinuous points
from Proposition \ref{Enul} we have $I^-_{\mu_2}=0$.
 From Proposition \ref{Enul} and
considering that $S(\mu_2 )=X_2$  we can compute the value of
$\lambda^-_{\mu_2}$ if we find a point $x$ such that
$\limsup_{n\to\infty}\frac{1}{n}\Lambda^-_n(x)$ be maximum. Denote by
$I^{-F}_n$ the map $I^-_n$ associated with the automaton $F$ and $I^{-F_1}_n$,
$I^{-F_2}_n$ those associated respectively with $F_1$ and $F_2$. Similarly
$\Lambda_n^{-F}$, $\Lambda_n^{-F_1}$ and $\Lambda_n^{-F_2}$ are the maps
$\Lambda^-_n$ associated with $F$, $F_1$ and $F_2$. As $F$ is the product of
$F_1$ by $F_2$ we have 
\small
$$
I^{-F}_n(x)=\max\{I^{-F_1}_n(x_1),I^{-F_2}_n(x_2)\}\mbox{ and
}\Lambda_n^{-F}(x)= \max\{\Lambda_n^{-F_1}(x_1),\Lambda_n^{-F_2}(x_2)\}.
$$
\normalsize
Remembering that $I_{n,\mu}^{-F}=\int_X I^{-F}_n(x)d\mu_1(x_1) d\mu_2(x_2)\,
$ then 

\[
 I_{n,\mu}^{-F}=\int_X
\max\{I^{-F_1}_n(x_1),I^{-F_2}_n(x_2)\}d\mu_1(x_1) d\mu_2(x_2).
\]
 If we
consider $F_2$ as a map on $X$ we can say that for $\mu_1\mu_2$-almost all $x\in
X$, $I^{-F_2}_n(x_1)$ is bounded. It follows that $\liminf\frac{1}{n}\int_X
I_n^{-F_2}(x) d\mu_1d\mu_2(x)=0, $
 hence 
\[
 I_\mu^-=\liminf\frac{1}{n}\int_X
I_n^{-F_1}(x)d\mu_1d\mu_2(x) =I_{\mu_1}^-=1.
\]
 If we consider successively
$F_1$ and $F_2$ as maps on $X$, we can see that for $\mu_1\mu_2$-almost all
$x\in
X$ we have $\Lambda_n^{-F_1}(x)=n$ and $\Lambda_n^{-F_2}(x)=rn$ which implies
that 
$$
\lambda_\mu^-=\liminf\int_X \max\{\Lambda^{-F_1}_n,\Lambda^{-F_2}_n\}
d\mu_1d\mu_2 =\liminf\int_X\Lambda^{-F_1}_n d\mu_1d\mu_2=I^-_{\mu_2}=r.
$$
Denoting by $\sigma_2$ the shift on $X_2$, considering that $h_\mu (F_1)=\log 2$
and $h_\mu (\sigma_2)=\log 3$ then $h_\mu (\sigma )=\log 2+\log 3$.
 From
Proposition \ref{Enul} and taking in account that  $F_2$ has equicontinuous
points we get $h_{\mu_2}(F_2)=0$. Considering successively $F_1$ and $F_2$ as
automata on their respective configuration spaces and on $X$ we can assert that
\[
h_\mu (F)=h_{\mu}(F_1)+h_{\mu}(F_2)=h_{\mu_1} (F_1)+h_{\mu_2} (F_2) =\log2.
\]
Finally applying inequalities (1) and (2) to $F$, we can conclude that for
this example the average Lyapunov exponents give better bounds for the entropy.
From Corollary \ref{corfi} one has 
$ h_\mu (F)=\log 2\le h_\mu (\sigma )(\lambda^+_\mu
+\lambda^-_\mu)=(\log 2+\log 3)r
$
\\ and from Theorem \ref{principal} we get 
$ h_\mu
(F)=\log 2\le h_\mu (\sigma )(I^+_\mu +I^-_\mu )=(\log 2+\log 3).
$

 In this
example $\mu$ is the uniform measure on $X$ so from Corollary \ref{corfi} the real $ h_\mu
(\sigma )(\lambda^+_\mu +\lambda^-_\mu)=(\log 2+\log 3)r$ is an upper bound of
the topological entropy of $F$. The topological entropy of $F$ is the sum of the
entropy of $F_1$ and $F_2$ and is equal to $(r+1)\log (2)$ which means that in
this case the topological inequality is strict.
\begin{premark}
It will be interesting to find no trivial examples for which inequality (2)
becomes an equality and with a strict inequality (1). The Proposition
\ref{Enul} suggest that we know very
little about sensitive cellular automata. The condition $h_\mu (\sigma
)(I^+_\mu +I^-_\mu )>0$ does not imply that $h_\mu (F)>0$.
\end{premark}
%

I am indebted to Fran\c{c}ois Blanchard for many stimulating conversations and
for his help for the writing.
An important part of this work has been done in the University of Chile in
Santiago, Mathematics's laboratory of the section Civil Engineering. I want
to thanks "FONDAP en Matematicas
Aplicadas, proyecto Modelamiento Estocastico" and ECOS for
the financial
support and Alejandro Maass for his numerous suggestions.

\end{document}